\documentclass[a4paper,12pt]{amsart}
\usepackage{amssymb}

\numberwithin{equation}{section}
\newtheorem{thm}{Theorem}[section]

\theoremstyle{remark}

\author[J.~Benameur]{Jamel Benameur}
\address{Department of Mathematics, College of Science, King Saud University\\
Riyadh 11451, Kingdom of Saudi Arabia}
\email{\sl jbenameur@ksu.edu.sa}
\author[M.~Benhamed]{Moez Benhamed}
\address{University of Tunis El Manar, Faculty of Sciences of Tunis, LR03ES04, 2092 Tunis, Tunisia}
\email{\sl moez.benhamed@hotmail.com}

\title[Long time decay to the solution to the 2D-DQG Equation]
{Long time decay to the solution to the 2D-DQG Equation}
\date{\today}

\begin{document}
\begin{abstract} In \cite{JAMH1}, we prove the well posedness of the
quasi-geostrophic equation $(QG)_{\alpha}\;,1/2<\alpha\leq 1$, in the space introduced by Z. Lei and F. Lin in \cite{ZY1}. In this
chapter we discuss the long time behaviour. Mainly, we prove that $\|\theta\|_{{\mathcal X}^{1-2\alpha}}$ decays to zero as time goes to infinity.
\end{abstract}


\subjclass[2000]{35A35; 35B05; 35B40; 35Q55; 81Q20}
\keywords{surface quasi-geostrophic equations; sub-critical case; critical spaces; long-time behavior}

\maketitle


\section{Introduction}
\label{S1}
In this paper, we study the initial value-problem for the two-dimensional
quasi-geostrophic equations with sub-critical dissipation $(QG)_{\alpha}$,
$$\left\{\begin{array}{l}
  \displaystyle\partial_t
\theta +u.\nabla \theta +k\Lambda ^{2\alpha}\theta=0, \\
u = R^{\bot}\theta=(-R_{2}\theta,R_{1}\theta), \\
\theta(0, .)=\theta_0,
\end{array}\right.
\leqno(QG)_{\alpha}$$
Part1Here $1/2<\alpha\leq 1$  is a real number and $k > 0$ is a dissipative coefficient. $\Lambda$ is the Riesz potential operator defined by the fractional power of $-\Delta$ $$\Lambda=(-\Delta)^{1/2}\quad \mathrm{and}\quad \widehat{\Lambda g}=\widehat{(-\Delta)^{1/2}g}=|\xi|\widehat{g}.$$ Thus,
 $$\widehat{\Lambda^{2\alpha}g}=\widehat{(-\Delta)^{\alpha}g}=|\xi|^{2\alpha}\widehat{g},$$ where $\widehat{g}(\xi)$ denotes the Fourier transform of $g$.\\
$\theta(x, t)$ is an unknown scalar function representing potential temperature. $u(x, t)$ is the velocity field determined by the scalar
stream function $\psi$ through
\begin{equation}u=(u_1, u_2)=\Big(-\frac{\partial}{\partial x_2}\psi, \frac{\partial}{\partial x_1}\psi\Big)=\nabla^{\bot}\psi.\label{e111}\end{equation} And the temperature $\theta$ and the stream function $\psi$ are related by $$\Lambda \psi=-\theta,$$ the relations in \eqref{e111} can be also combined into
\begin{eqnarray*} u=-\nabla^{\bot}(-\Delta)^{-1/2}\theta=- R^{\bot}\theta=(\partial_2(-\Delta)^{-1/2}\theta,\; -\partial_1(-\Delta)^{-1/2}\theta, \end{eqnarray*}
where $R_j,\, j = 1, 2$ is the $2D$ Riesz transform defined by $$\widehat{R_{j}f}=-\frac{\xi_j}{|\xi|}\widehat{f}.$$

The $2D$ quasi-geostrophic fluid is an important model in geophysical fluid dynamics, they are special cases of the general
quasi-geostrophic approximations for atmospheric and oceanic fluid flow with the small local Rossby number which ensures
the validity of the geostrophic balance between the pressure gradient and the Coriolis force (see Pedlosky \cite{JP}). Furthermore,
this quasi-geostrophic fluid motion equation shares many features with fundamental fluid motion equations. When $k = 0$,
this equation is comparable to the vorticity formulation of the Euler equations (see Majda and Tabak \cite{CMT}), and $(QG)_{\alpha}$ with
$\alpha = 0$ is similar to a non viscous wind driven circulation equation (see Pedlosky \cite{JP}). What is more, $(QG)_{\alpha}$ with $\alpha = 1$ shares
similar features with the three-dimensional Navier-Stokes equations. Thus $\alpha = 1$ is therefore referred as the critical case,
while the cases $\alpha >1,\, \alpha = 1 \;\mathrm{and}\; \alpha <1$  are supercritical and subcritical, respectively.\\

However, it is desirable to understand the  global existence for the dissipative quasi-geostrophic
equation. Existence of a global weak solution was established by Resnick \cite{SR}. Furthermore, in the subcritical case, Constantin
and Wu \cite{CW11} proved that every sufficiently smooth initial data give rise to a unique global smooth solution. In the
critical case, $\alpha = 1/2$, Constantin, Cordoba and Wu \cite{CWDC} established the existence of a unique global classical solution
corresponding to any initial data that are small in $L^\infty$. The assumption requiring smallness in $L^\infty$ was removed by
Caffarelli and Vasseur \cite{LACAV}, Kiselev, Nazarov and Volberg \cite{KNV} and Dong and Du \cite{HDDD}. In \cite{KNV} the authors proved
persistence of a global solution in $\mathcal{C}^\infty$ corresponding to any $\mathcal{C}^\infty$ periodic initial data. Chae and
Lee \cite{CDLJ} obtained the global existence and uniqueness of solutions for data in the Besov space $B^{2-\gamma}_
{2, 1}$ with a small $\dot{B}^{2-\gamma}_
{2, 1}$ norm.\\

Recently,  Benameur and  Benhamed \cite{JAMH1} proved the global existence for the dissipative quasi-geostrophic for a small initial data. They introduce the critical space defined by \begin{equation}\label{X} {\mathcal X}^{\sigma}(\mathbb R^2)=\Big\{f\in \mathcal{S'}(\mathbb R^2);\;(\xi\mapsto|\xi|^{\sigma}\widehat{f}(\xi))\in L^1(\mathbb R^2) \Big\}.\end{equation}
More precisely
\begin{thm}\label{thjm1} Let $\theta^0\in{\bf{\mathcal X}}^{1-2\alpha}(\mathbb R^2)$. There is a time $T>0$ and unique solution
$\theta\in\mathcal C([0,T],{\bf{\mathcal X}}^{1-2\alpha}(\mathbb R^2))$ of $(QG)_{\alpha}$, moreover $\theta\in L^1([0,T],{\bf{\mathcal X}}^{1}(\mathbb
R^2))$. If $\|\theta^0\|_{{\bf{\mathcal X}}^{1-2\alpha}}<1/4$, the solution is global and
\begin{equation}\label{thjm1eq}
\|\theta\|_{{\bf{\mathcal X}}^{1-2\alpha}}+\frac{1-4\|\theta^0\|_{{\bf{\mathcal X}}^{1-2\alpha}}}{2}\int_0^t\|\theta\|_{{\bf{\mathcal X}}^{1}}\leq
\|\theta^0\|_{{\bf{\mathcal X}}^{1-2\alpha}},\;\;\forall t\geq0.
\end{equation}
\end{thm}
The main purpose of this work is study the long time limit of the solutions of two-dimensional
quasi-geostrophic equation $(QG)_{\alpha}$ when $1/2<\alpha<1$, in this case Niche and Schonbek \cite{CNES} prove that if the initial data $\theta^0$ is in
$L^2(\mathbb R^2)$, then the $L^2$ norm of the solution tends to zero but with no uniform rate, that is, there are
solutions with arbitrary slow decay. If $\theta^0\in L^p(\mathbb R^2)$,  with $1\leq  p\leq 2$, they obtain a uniform
decay rate in $L^2$.\\

 We state now our main result.
\begin{thm}\label{thjm2}
Let $2/3<\alpha <1$ and $\theta\in\mathcal C(\mathbb R^+,{\bf{\mathcal X}}^{1-2\alpha}(\mathbb R^2))$ be a global solution of $(QG)_{\alpha}$ given by
Theorem \ref{thjm1}. Then
\begin{equation}\label{thjm2eq}
\displaystyle \lim_{t\rightarrow \infty}\|\theta\|_{{\bf{\mathcal X}}^{1-2\alpha}}=0.
\end{equation}
\end{thm}

\end{document}